\documentclass[10pt,a4paper]{article}
\usepackage{amsfonts}
\usepackage{}
\usepackage{amsmath}
\usepackage{latexsym}
\usepackage{amssymb}
\usepackage{amstext}
\usepackage{fancyhdr}
\usepackage{indentfirst}
\usepackage{graphicx}
\usepackage{epsfig}
\usepackage{pgf}
\usepackage{listings}
\usepackage{cases}
\usepackage{epstopdf}
\usepackage{graphicx}
\usepackage{enumerate}
\date{}

\begin{document}
\title{Quantum State Transfer on Neighborhood Corona of Two Graphs}
\author{Xiao-Qin Zhang$^a$, Qi Xiong$^a$, Gui-Xian Tian$^{a}$\footnote{Corresponding author. E-mail address: gxtian@zjnu.cn (G-X. Tian)}, Shu-Yu Cui$^b$\\\\
    {\small{\it $^a$Department of Mathematics,}}
    {\small{\it Zhejiang Normal University, Jinhua, 321004, China}}\\
    {\small{\it $^b$Xingzhi College, Zhejiang Normal University, Jinhua, 321004, China}}
}\maketitle
\begin{abstract}
Given two graphs $G_{1}$ of order $n_{1}$ and $G_{2}$, the neighborhood corona of $G_{1}$ and $G_{2}$, denoted by $G_{1}\bigstar G_{2}$, is the graph obtained by taking one copy of $G_{1}$ and taking $n_{1}$ copies of $G_{2}$, in the meanwhile, linking all the neighbors of the $i$-th vertex of $G_{1}$ with all vertices of the $i$-th copy of $G_{2}$. In our work, we give some conditions that $G_{1}\bigstar G_{2}$ is not periodic. Furthermore, we demonstrate some sufficient conditions for $G_{1}\bigstar G_{2}$ having no perfect state transfer. Some examples are provided to explain our results. In addition, for the reason that the graph admitting perfect state transfer is rare, we also consider pretty good state transfer on neighborhood corona of two graphs. We show some sufficient conditions for $G_{1}\bigstar G_{2}$ admitting pretty good state transfer.

\emph{AMS classification:} 05C50 81P68

\emph{Keywords:} Quantum state transfer; Neighborhood Corona;
Spectrum; Perfect state transfer; Pretty good state transfer
\end{abstract}

\section{Introduction}

Quantum walk is a natural generalization of classical random walk on graphs. In $1998$, Farhi and Gutmann \cite{Farhi} firstly put forward the concept of continuous-time quantum walk. Given a graph $G$, let $A(G)$, $L(G)$ and $Q(G)$ be its adjacency matrix, Laplacian matrix and signless Laplacian matrix, respectively. Suppose that $X(G)$ is a Hermitian matrix associated with $G$, then the unitary matrix $H_{G}(t)=\text{exp}(-itX(G))$ is the transition matrix of continuous-time quantum walk corresponding to $X(G)$ for $i^2=-1$ and $t>0$, where $X(G)$ may be $A(G)$, $L(G)$, $Q(G)$ and so on. In 2003, Bose \cite{Bose} studied the task of information transition in a quantum spin system. Christandl et al. \cite{Christandl} showed that this task can be lessened to the question of perfect state transfer. Let $e_{x}^{n}$ be a column vector of order $n$ whose element of $x$-th position is $1$ and $0$, otherwise. Sometimes, $e_{x}^{n}$ can be recorded briefly by $e_{x}$.
If
\begin{equation}
\text{exp}(-itX(G))e_{x}=\lambda e_{y}
\end{equation}
with $|\lambda| =1$ for two vertices $x, y$ of $G$, then we say that $G$ admits \textit{perfect state transfer} (PST for short) relative to matrix $X(G)$ between $x$ and $y$ at time $t$. Particularly, if $X(G)=A(G)$ (resp. $L(G)$ or $Q(G)$), then $G$ admits \textit{perfect state transfer} (resp. \textit{Laplacian perfect state transfer} (LPST for short) or \textit{signless Laplacian perfect state transfer} (SLPST for short)).

Since the graph admitting PST is rare. Godsil \cite{Godsil2012} proposed a new concept called pretty good state transfer (PGST for short), whose restriction is more relaxing than that of PST. The transpose of $x$ is denoted as $x^{T}$.
if
\begin{equation}
|e_{y}^{T}\text{exp}(-i\tau X(G))e_{x}|>1-\epsilon,
\end{equation}
for any $\epsilon>0$, then we say that $G$ admits \emph{pretty good state transfer} (PGST for short) between two vertices $x$ and $y$ at time $\tau$. Similarly, if $X(G)=A(G)$ (resp. $L(G)$ or $Q(G)$), then $G$ admits \textit{pretty good state transfer} (resp. \textit{Laplacian pretty good state transfer}(LPGST for short) or \textit{signless Laplacian pretty good state transfer} (SLPGST for short)).

Recently, many articles focus on PST and PGST of composite graphs which are obtained by some graph operations. For example, Li et al. \cite{Yipeng Li} considered LPST and LPGST of $\mathcal{Q}$-graph. They showed that, for an $r$-regular graph $G$, if $r+1$ is a prime, then the $\mathcal{Q}$-graph of $G$ has no LPST but admits LPGST. Recently, Zhang et al. \cite{zhang} studied SLPST and SLPGST in $\mathcal{Q}$-graph. Ackelsberg et al. \cite{laplaciancorona} considered LPST and LPGST of corona graphs. They showed that corona graph $G\circ H$ has no LPST, but it occurs LPGST under some special conditions. In $2017$, Ackelsberg et al. \cite{adjacentcorona} considered PST and PGST of corona graphs. They showed that $G\circ K_{n}$ has no PST and $G\circ K_{1}$ admits PGST under some suitable conditions. In $2021$, Tian et al. \cite{Tian2021} considered SLPST and SLPGST of corona graphs. They demonstrated that $K_{2}\circ H$ has no SLPST between the two vertices of $K_{2}$. They also showed that $G \circ \overline{K_{m}}$ admits SLPGST under some suitable condition. In $2021$, Wang and Liu \cite{edge complemented coronas} considered LPST and LPGST of edge complemented coronas. They gave some sufficient conditions such that edge complemented corona $G \diamond H$ has no LPST. They also showed that $G \diamond H$ admits LPGST under some suitable conditions. Recently, Li et al. \cite{neighborhood coronas} gave some sufficient conditions for extended neighborhood coronas to
have Laplacian perfect state transfer. For more details about these directions, readers may refer to \cite{laplacian,jiweijiuhui,Coutinho2015,Fan,Xiaogang Liu,zhangliu}.

Motivated by the aforementioned results, we mainly focus on PST and PGST of neighborhood corona of two graphs. Given two graphs $G_{1}$ of order $n_1$ and $G_{2}$, the \emph{neighborhood corona} of $G_{1}$ and $G_{2}$, denoted by $G_{1}\bigstar G_{2}$, is the graph obtained by taking one copy of $G_{1}$ and $n_1$ copies of $G_{2}$, in the meanwhile, linking all the neighbors of the $i$-th vertex of $G_{1}$ with all vertices of the $i$-th copy of $G_{2}$. In our work, we first give some conditions that $G_{1}\bigstar G_{2}$ is not periodic. With the help of these conditions, we demonstrate some sufficient conditions for $G_{1}\bigstar G_{2}$ having no PST. Furthermore, some examples are provided to explain our results. Finally, for the reason that the graph admitting PST is rare, we also consider PGST on neighborhood corona of two graphs. We show some sufficient conditions for $G_{1}\bigstar G_{2}$ admitting PGST. It turns out that $G_{1}\bigstar \overline{K_{n_{2}}}$ with $G_{1}$ admitting PST, ${C_{4}} \bigstar \overline{K_{n_{2}}}$ and $C_{4}\bigstar G$ have PGST under some special conditions.

\section{Preliminaries}
Throughout this paper, we only consider undirected simple graphs. Let $j_{n}$ and $J_{n}$ be the all-one column vector of order $n$ and all-one square matrix of order $n$, respectively. Let $[n]$ be the set of $\{1, 2, \ldots, n\}$.

Suppose that $G$ is a graph of order $n$ and $\lambda_{j}$ is the eigenvalue of $A(G)$ with multiplicity $l_{j}$ for $j\in[p]$, where $l_1+\cdots+l_p=n$. The spectrum of $A(G)$ is denoted by $\text{Sp}(G)$, then $\text{Sp}(G)=\{\lambda_{j}^{l_{j}}:j\in[p]\}$.
Let $\{{x_{1}^{(j)}}, x_{2}^{(j)}, \ldots,x_{l_{j}}^{(j)}\}$ be an orthonomal basis of the eigenvalue space $V_{\lambda_{j}}$ of $\lambda_{j}$. The eigenprojector of $\lambda_{j}$, denoted by $f_{\lambda_{j}}(G)$, is $f_{\lambda_{j}}(G)=\sum\nolimits_{i = 1}^{l_j}x_{i}^{(j)}(x_{i}^{(j)})^{T}$ and $\sum\nolimits_{j=1}^{p}f_{\lambda_{j}}(G)=I$. Obviously, $f_{\lambda_{j}}(G)f_{\lambda_{k}}(G)=0$ for $j\neq k$ and $f_{\lambda_{j}}(G)$ is an idempotent matrix.
According to the eigenprojectors, we get the spectral decomposition of $A(G)$, i.e.
\begin{equation}
A(G)=A(G)\sum\limits_{j=1}^{p}f_{\lambda_{j}}(G)=\sum\limits_{j=1}^{p}\lambda_{j}f_{\lambda_{j}}(G).
\end{equation}
Hence,
\begin{equation}
H_{G}(t)=\text{exp}(-itA(G))=\sum\limits_{j=1}^{p}e^{-it\lambda_{j}}f_{\lambda_{j}}(G).
\end{equation}
The eigenvalue support of vertex $x$, denoted by $\text{supp}_{G}(x)$, is the set of eigenvalues $\lambda_{j}$ such that $f_{\lambda_{j}}(G)e_{x}\neq0$. Two vertices $x$ and $y$ are \textit{strongly cospectral} whenever $f_{\lambda_{j}}(G)e_{x}=\pm f_{\lambda_{j}}(G)e_{y}$ for any $j\in[p]$.
Let $S^{+}=\{\lambda:f_{\lambda}(G)e_{x}=f_{\lambda}(G)e_{y}\}$ and $S^{-}=\{\lambda:f_{\lambda}(G)e_{x}=-f_{\lambda}(G)e_{y}\}$.

\paragraph{Theorem 2.1.}(Coutinho \cite{Coutinho2014}) Assume that $G$ is a graph with vertex set satisfying $|V(G)|\geq 2$, and $u, v\in V(G)$.
If $\lambda_{0}$ is the maximum eigenvalue of $G$,
then $G$ admits PST between the vertices $u$ and $v$ if and only if the following conditions hold.
\begin{enumerate}[(i)]
 \item Two vertices, $u$ and $v$, are strongly cospectral.
 \item Non-zero elements in $\text{supp}_{G}(u)$ are either all integers or all quadratic integers. Moreover, for each eigenvalue $\lambda\in \text{supp}_{G}(u)$, there exists a square-free integer $\Delta$ and integers $a$, $b_\lambda$ such that,
  \begin{equation*}
  \lambda=\frac{1}{2}(a+b_\lambda\sqrt{\Delta}).
  \end{equation*}
  Here we allow $\Delta=1$ if all eigenvalues in $\text{supp}_{G}(u)$ are integers, and $a=0$ if all eigenvalues in $\text{supp}_{G}(u)$ are multiples of $\sqrt{\Delta}$.
 \item $\lambda\in{S^+}$ if and only if $\dfrac{\lambda_0-\lambda}{g\sqrt{\Delta}}$ is even and  $\lambda\in{S^-}$ if and only if $\dfrac{\lambda_0-\lambda}{g\sqrt{\Delta}}$ is odd, where
 $$g=\gcd\left(\left\{\dfrac{\lambda_0-\lambda}{\sqrt{\Delta}}:\lambda\in \text{supp}_{G}(u)\right\}\right).$$
\end{enumerate}
Moreover, if the conditions above hold, then the following also hold.
\begin{enumerate}[(1)]
     \item There exists a minimum time $\tau_0>0$ at which PST occurs between $u$ and $v$, and
     \begin{equation*}
     \tau_0=\frac{1}{g}\dfrac{\pi}{\sqrt{\Delta}}.
     \end{equation*}
     \item The time of PST, $\tau$ is an odd multiple of $\tau_0$.
     \item The phase of PST is given by $\lambda=e^{-it\lambda_0}$.
\end{enumerate}

In order to characterize graphs admitting PST (or PGST), the following two lemmas play a crucial role in our study process.
\paragraph{Lemma 2.2.}(Godsil \cite{Godsil}) If a graph $G$ admits PST between two vertices $u$ and $v$ at time $t$, then $G$ is periodic at vertex $u$ (or $v$) at time $2t$.

\paragraph{Lemma 2.3.}(Godsil \cite{equation}) A graph $G$ at vertex $v$ is periodic if and only if one of the following conditions holds:
\begin{enumerate}[(i)]
\item all elements of $\text{supp}_{G}(v)$ are integers;
\item for each eigenvalue of $\text{supp}_{G}(v)$, there is a square-free integer $\Delta$, integer $a$ and corresponding some integer $b_{\lambda}$ so that $\lambda=\frac{1}{2}(a+b_{\lambda}\sqrt{\Delta})$.
\end{enumerate}

\paragraph{Theorem 2.4.}(Hardy and Wright \cite{numbertheory}) Assume that $1, \lambda_1, \ldots, \lambda_m$ are linearly independent over $\mathbb{Q}$. Then, for any real numbers $\alpha_1, \ldots, \alpha_m$ and $N>0$, $\epsilon>0$, there exist integers $\alpha>N$ and $\gamma_1, \ldots, \gamma_m$ such that
\begin{equation}\label{4*}
\vert{\alpha\lambda_{k}-\gamma_k-\alpha_k}\vert<\epsilon,
\end{equation}
for each $k\in[m]$. Equivalently, (\ref{4*}) can be restated by $\alpha\lambda_{k}-\gamma_k\approx \alpha_k$ for omitting the dependence on $\epsilon$.

\paragraph{Lemma 2.5.}(Richards \cite{galoistheory}) The set $\{\sqrt{\Delta}:\Delta$ is a square-free integer$\}$ is linearly independent over
$\mathbb{Q}$.

\paragraph{Lemma 2.6.}(Coutinho \cite{Coutinho2014}) A real number
$\lambda$ is a quadratic integer if and only if there exist integers $a$, $b$ and $\Delta$ such that $\Delta$ is square-free and one of the following cases holds:
\begin{enumerate}[(i)]
    \item $\lambda=a+b\sqrt{\Delta}$ and $\Delta\equiv2,3\;(\text{mod}\;4)$;
    \item $\lambda=\frac{1}{2}(a+b\sqrt{\Delta}),\; \Delta\equiv1\;(\text{mod}\;4)$, and $a$ and $b$ have the same parity.
\end{enumerate}

\section{Neighborhood Corona of Graphs}
We denote the neighborhood corona of graph $G_{1}$ and $G_{2}$ by $G_{1} \bigstar G_{2}$. Let binary ordered pair $V(G_{1} \bigstar G_{2})= V(G_{1})\times ({0}\cup V(G_{2}))$ be the vertex set of $G_{1} \bigstar G_{2}$ and $|V(G_{1})|=n_{1}$, $|V(G_{2})|=n_{2}$. According to the definition of $G_{1} \bigstar G_{2}$, the adjacency relation is given by:
\begin{equation}
((x,y),(x',y'))\in{E(G_{1} \bigstar G_{2})}\Longleftrightarrow
{\Bigg\lbrace{
\begin{aligned}
&y=y'=0 \quad \text{and} \quad (x,x')\in{E(G_{1})} \quad\quad\; \text{or}\\
&x=x'  \quad \text{and} \quad (y,y')\in{E(G_{2})} \quad\quad\quad\quad \text{or}\\
&(x,x')\in{E(G_{1})} \quad  \text{and exactly one of}\; y \; \text{and} \; y'\;\text{ is} \; 0.
\end{aligned}}
}
\end{equation}
According to the adjacency relation above, we easily get the adjacency matrix of $G_{1} \bigstar G_{2}$ as follows:
\begin{equation}\label{1}
A(G_{1}\bigstar G_{2}) = \left( {\begin{array}{*{20}{c}}
{A(G_{1})}&{A(G_{1}) \otimes j_{{n_2}}^T}\\
{A(G_{1}) \otimes {j_{{n_2}}}}&{{I_{{n_1}}} \otimes A(G_{2})}
\end{array}} \right),
\end{equation}
where $``\otimes"$ denotes the Kronecker product of two matrices.

At first, we recall the adjacency spectrum of $G_{1}\bigstar G_{2}$, which is attributed to Gopalapillai in \cite{neighborhood}.\\
\\
\textbf{Theorem 3.1.} (Gopalapillai \cite{neighborhood}) Let $G_{1}$ be any graph of order $n_{1}$ and $G_{2}$ be any $k$-regular connected graph with $n_{2}$ vertices. Let $\text{Sp}(G_{1})=\{\lambda_{i}^{l_{i}} | i\in [p]\}$ and $\text{Sp}(G_{2})=\{\eta_{j}^{l_{j}'} | j\in [q]\}$, where the power represents the multiplicity. Then the eigenvalues of $A(G_{1}\bigstar G_{2})$ are
\begin{enumerate}[(1)]
\item
\begin{equation}\label{6}
\lambda_{i+}=\frac{\lambda_{i}+k +\sqrt{(\lambda_{i}-k)^{2}+4n_{2}\lambda_{i}^{2}}}{2}
\end{equation}
and
\begin{equation}\label{7}
\lambda_{i-}=\frac{\lambda_{i}+k -\sqrt{(\lambda_{i}-k)^{2}+4n_{2}\lambda_{i}^{2}}}{2}
\end{equation}
with multiplicity $l_{i}$ for $i\in[p]$.
\item
$\eta_{j}$ with multiplicity $n_{1}$ for $j\in [q]\setminus\{1\}$.
\end{enumerate}

In particular, if there is some $i_{0}\in [p]$ such that $\lambda_{i_{0}}=0$, then
    \begin{equation}
     \lambda_{i_{0}+}=\frac{\lambda_{i_{0}}+k +\sqrt{(\lambda_{i_{0}}-k)^{2}+4n_{2}\lambda_{i_{0}}^{2}}}{2}=k
    \end{equation}
    and
    \begin{equation}
     \lambda_{i_{0}-}=\frac{\lambda_{i_{0}}+k -\sqrt{(\lambda_{i_{0}}-k)^{2}+4n_{2}\lambda_{i_{0}}^{2}}}{2}=0
    \end{equation}
  with multiplicity $l_{i_{0}}$.\\
\\
\textbf{Proposition 3.2}. Let $G_{1}$ be any graph with $n_{1}$ vertices and $G_{2}$ be a $k$-regular graph with $n_{2}$ vertices. Let $\text{Sp}(G_{1})=\{{\lambda_{i}}^{l_{i}} |i\in [p]\}$ and $\text{Sp}(G_{2})=\{{\eta_{j}}^{l'_{j}} | j\in [q]\}$, where the respective multiplicities of $\lambda_{i}$ and $\eta_{j}$ are $l_{i}$ and $l'_{j}$ for $i\in [p]$, $j\in [q]$.
\begin{enumerate}[(i)]
   \item If $\lambda_{i}\neq 0$ for any $i\in[p]$, then
    \begin{enumerate}[(1)]
    \item
    the eigenprojector of eigenvalue $\lambda_{i\pm}$ of  $A(G_{1}\bigstar G_{2})$ is
    \begin{equation}\label{3}\small
    {F_{{\lambda _{i \pm }}}}(G_{1}\bigstar G_{2}) = \frac{{{\lambda _i}^2}}{{{{({\lambda _{i \pm }} - k)}^2} + {n_2}{\lambda _i}^2}}\left( {\begin{array}{*{20}{c}}
    {\frac{{{{({\lambda _{i \pm }} - k)}^2}}}{{{\lambda _i}^2}}{f_{{\lambda _i}}}({G_1})}&{\frac{{({\lambda _{i \pm }} - k)}}{{{\lambda _i}}}{f_{{\lambda _i}}}({G_1}) \otimes {j_{{n_2}}}^T}\\
    {\frac{{({\lambda _{i \pm }} - k)}}{{{\lambda _i}}}{f_{{\lambda _i}}}({G_1}) \otimes {j_{{n_2}}}}&{{f_{{\lambda _i}}}({G_1}) \otimes {J_{{n_2}}}}
     \end{array}} \right)
    \end{equation}
    for $i\in[p]$;
    \item
    the eigenprojector of eigenvalue $\eta_{j}$ of $A(G_{1}\bigstar G_{2})$ is
     \begin{equation}\label{5}
    {F_{{\eta _j}}}(G_{1}\bigstar G_{2}) = \left( {\begin{array}{*{20}{c}}
     0&0\\
     0&{{I_{{n_1}}} \otimes {(f_{{\eta _j}}}({G_2})-\delta_{\eta _j,k}\frac{1}{n_{2}}J_{n_{2}})}    \end{array}} \right)
     \end{equation}
    for $j\in[q]$.
    \end{enumerate}
    Hence the spectrum decomposition of $A(G_{1}\bigstar G_{2})$ is given by:
 \begin{equation}\label{2}
A({G_1} \bigstar {G_2}) = \sum\limits_{i = 1}^p {\sum\limits_ \pm  {{\lambda _{i \pm }}{F_{{\lambda _{i \pm }}}}({G_1} \bigstar {G_2})} }  + \sum\limits_{j = 1}^q {{\eta _j}{F_{{\eta _j}}}} ({G_1} \bigstar {G_2}).
\end{equation}
    \item If there is some $i_{0}\in [p]$ such that $\lambda_{i_{0}}=0$, then
    \begin{enumerate}[(1)]
    \item
    the eigenprojector of eigenvalue $\lambda_{i\pm}$ of $A(G_{1}\bigstar G_{2})$ is
    \begin{equation}\small
    {F_{{\lambda _{i \pm }}}}(G_{1}\bigstar G_{2}) = \frac{{{\lambda _i}^2}}{{{{({\lambda _{i \pm }} - k)}^2} + {n_2}{\lambda _i}^2}}\left( {\begin{array}{*{20}{c}}
    {\frac{{{{({\lambda _{i \pm }} - k)}^2}}}{{{\lambda _i}^2}}{f_{{\lambda _i}}}({G_1})}&{\frac{{({\lambda _{i \pm }} - k)}}{{{\lambda _i}}}{f_{{\lambda _i}}}({G_1}) \otimes {j_{{n_2}}}^T}\\
    {\frac{{({\lambda _{i \pm }} - k)}}{{{\lambda _i}}}{f_{{\lambda _i}}}({G_1}) \otimes {j_{{n_2}}}}&{{f_{{\lambda _i}}}({G_1}) \otimes {J_{{n_2}}}}
     \end{array}} \right)
    \end{equation}
    for $i\in[p]$ and $i\neq i_{0}$;
    \item
     the eigenprojector of eigenvalue $\eta_{j}$ of $A(G_{1}\bigstar G_{2})$ is
     \begin{equation}
    {F_{{\eta _j}}}(G_{1}\bigstar G_{2}) = \left( {\begin{array}{*{20}{c}}
     0&0\\
     0&{{I_{{n_1}}} \otimes {f_{{\eta _j}}}({G_2})}
    \end{array}} \right)
     \end{equation}
    for $j\in[q]$ and $j\neq1$;
    \item
    the eigenprojector of eigenvalue $k$ of $A(G_{1}\bigstar G_{2})$ is
    \begin{equation}
    {F_k(G_{1}\bigstar G_{2})} = \frac{1}{n_{2}}\left( {\begin{array}{*{20}{c}}
    0&0\\
   0&{{f_0}({G_1}) \otimes {J_{{n_2}}}}
    \end{array}} \right)+\left( {\begin{array}{*{20}{c}}
     0&0\\
     0&{{I_{{n_1}}} \otimes {(E_{{k}}}({G_2})-\frac{1}{n_{2}}J_{n_{2}})}
    \end{array}} \right).
   \end{equation}
   \item
   the eigenprojector of eigenvalue $0$ of $A(G_{1}\bigstar G_{2})$ is
   \begin{equation}\label{9}
   {F_0(G_{1}\bigstar G_{2})} = \left( {\begin{array}{*{20}{c}}
{{f_0}({G_1})}&0\\
0&0
\end{array}} \right).
   \end{equation}
    \end{enumerate}
  Hence the spectrum decomposition of $A(G_{1}\bigstar G_{2})$ is given by:
  \begin{equation}\small
A({G_1} \bigstar {G_2}) = \sum\limits_{i \ne {i_0}}^{} {\sum\limits_ \pm  {{\lambda _{i \pm }}{F_{{\lambda _{i \pm }}}}({G_1} \bigstar {G_2})} }  + \sum\limits_{j = 2}^q {{\eta _j}{F_{{\eta _j}}}} ({G_1} \bigstar {G_2}) + k{F_k}({G_1} \bigstar {G_2}) + 0{F_0}({G_1} \bigstar {G_2}).
   \end{equation}
\end{enumerate}
\textit{Proof}: (i) Suppose that $\{x_{j}^{(i)}  | j\in[l_{i}]\}$ is the set of all orthonormal eigenvectors of eigenvalue $\lambda_{i}$ of $G_{1}$ and $\{y_{i}^{(j)}  | i\in[l_{j}']\}$ is the set of all orthonormal eigenvectors of eigenvalue $\eta_{j}$ of $G_{2}$ for $i\in [p]$, $j\in [q]$.
The eigenvectors of $A(G_{1}\bigstar G_{2})$ can be easily obtained from the proof of Theorem 2.1 in \cite{neighborhood}. For the convenience of readers, we give the detailed proof.
Let \[X_{{i\pm }}^j = \frac{{|{\lambda _i}|}}{{\sqrt {{{({\lambda _{{i\pm }}} - k)}^2} + {n_2}\lambda _i^2} }}\left( {\begin{array}{*{20}{c}}
{\frac{{{\lambda _{i \pm }} - k}}{{{\lambda _i}}}x_j^{(i)}}\\
{x_j^{(i)} \otimes {j_{{n_2}}}}
\end{array}} \right)\]  for $j\in [l_{i}]$, $i\in [p]$.
According to (\ref{1}), then
\[\begin{array}{l}
A({G_1} \bigstar {G_2})X_{i \pm }^j
 = \frac{{|{\lambda _i}|}}{{\sqrt {{{({\lambda _{i \pm }} - k)}^2} + {n_2}\lambda _i^2} }}\left( {\begin{array}{*{20}{c}}
{A({G_1})}&{A({G_1}) \otimes j_{{n_2}}^T}\\
{A({G_1}) \otimes {j_{{n_2}}}}&{{I_{{n_1}}} \otimes A({G_2})}
\end{array}} \right)\left( {\begin{array}{*{20}{c}}
{\frac{{{\lambda _{i \pm }} - k}}{{{\lambda _i}}}x_j^{(i)}}\\
{x_j^{(i)} \otimes {j_{{n_2}}}}
\end{array}} \right)\\
\;\;\;\;\;\;\;\; \;\;\;\;\;\;\;\;\;\;\;\;\;\;\;\;\;= \frac{{|{\lambda _i}|}}{{\sqrt {{{({\lambda _{i \pm }} - k)}^2} + {n_2}\lambda _i^2} }}\left( {\begin{array}{*{20}{c}}
{({\lambda _{i \pm }} - k + {n_2}{\lambda _i})x_j^{(i)}}\\
{{\lambda _{i \pm }}x_j^{(i)} \otimes {j_{{n_2}}}}
\end{array}} \right)\\
 \;\;\;\;\;\;\;\; \;\;\;\;\;\;\;\;\;\;\;\;\;\;\;\;\;= {\lambda _{i \pm }}\frac{{|{\lambda _i}|}}{{\sqrt {{{({\lambda _{i \pm }} - k)}^2} + {n_2}\lambda _i^2} }}\left( {\begin{array}{*{20}{c}}
{\frac{{{\lambda _{i \pm }} - k}}{{{\lambda _i}}}x_j^{(i)}}\\
{x_j^{(i)} \otimes {j_{{n_2}}}}
\end{array}} \right)\\
\;\;\;\;\;\;\;\; \;\;\;\;\;\;\;\;\;\;\;\;\;\;\;\;\;={\lambda _{i \pm }}X_{i \pm }^j
\end{array}\]
for $i\in [p]$ and $j\in [l_{i}]$, where the last equality holds because $(\lambda_{i\pm}-k+n_{2}\lambda_{i})\lambda_{i}=(\lambda_{i\pm}-k)\lambda_{i\pm}$.
Hence, according to the definition of eigenprojector, one gets
\begin{equation}\small
\begin{split}
{F_{{\lambda _{i \pm }}}}({G_1} \bigstar {G_2})&= \sum\limits_{j = 1}^{{l_i}} {X_{i \pm }^j{{(X_{i \pm }^j)}^T}}\\
&= \sum\limits_{j = 1}^{{l_i}} {\frac{{{\lambda _i}^2}}{{{{({\lambda _{i \pm }} - k)}^2} + {n_2}{\lambda _i}^2}}}
\left( {\begin{array}{*{20}{c}}
{\frac{{{\lambda _{i \pm }} - k}}{{{\lambda _i}}}x_j^{(i)}}\\
{x_j^{(i)} \otimes {j_{{n_2}}}}
\end{array}} \right){\left( {\begin{array}{*{20}{c}}
{\frac{{{\lambda _{i \pm }} - k}}{{{\lambda _i}}}x_j^{(i)}}\\
{x_j^{(i)} \otimes {j_{{n_2}}}}
\end{array}} \right)^T}\\
&= \frac{{{\lambda _i}^2}}{{{{({\lambda _{i \pm }} - k)}^2} + {n_2}{\lambda _i}^2}}
\left( {\begin{array}{*{20}{c}}
{\frac{{{{({\lambda _{i \pm }} - k)}^2}}}{{{\lambda _i}^2}}{f_{{\lambda _i}}}({G_1})}&{\frac{{({\lambda _{i \pm }} - k)}}{{{\lambda _i}}}{f_{{\lambda _i}}}({G_1}) \otimes {j_{{n_2}}}^T}\\
{\frac{{({\lambda _{i \pm }} - k)}}{{{\lambda _i}}}{f_{{\lambda _i}}}({G_1}) \otimes {j_{{n_2}}}}&{{f_{{\lambda _i}}}({G_1}) \otimes {J_{{n_2}}}}
\end{array}} \right).
\end{split}
\end{equation}
Let
$ Y_j^{ii'} = \left( {\begin{array}{*{20}{c}}
0\\
{{e_{i'}} \otimes y_i^{(j)}}
\end{array}} \right)$ for $ j\in [q]$ and $ i\in [l_{j}']$, $ i'\in [n_{1}]$ , where $e_{i'}$ is the characteristic vector of order $n_{1}$ and $y_i^{(j)}$ is a unit vector orthogonal to $j_{n_{2}}$.
Then
\[\begin{array}{l}
A({G_1} \bigstar {G_2})Y_j^{ii'}
 = \left( {\begin{array}{*{20}{c}}
{A({G_1})}&{A({G_1}) \otimes j_{{n_2}}^T}\\
{A({G_1}) \otimes {j_{{n_2}}}}&{{I_{{n_1}}} \otimes A({G_2})}
\end{array}} \right)\left( {\begin{array}{*{20}{c}}
0\\
{{e_{i'}} \otimes y_i^{(j)}}
\end{array}} \right)\\
 \;\;\;\;\;\;\;\; \;\;\;\;\;\; \;\;\;\;\;\;\;\; \;\;\;= {\eta _j}\left( {\begin{array}{*{20}{c}}
0\\
{{e_{i'}} \otimes y_i^{(j)}}
\end{array}} \right)\\
  \;\;\;\;\;\;\;\; \;\;\;\;\;\; \;\;\;\;\;\;\;\; \;\;\;= {\eta _j}Y_j^{ii'}.\\
\end{array}\]
Thus,
\[\begin{array}{l}
{F_{{\eta _j}}}({G_1} \bigstar {G_2})
 = \sum\limits_{i' = 1}^{{n_1}} {\sum\limits_{i = 1}^{{l_j}^\prime } {Y_j^{ii'}(Y_j^{ii'}} {)^T}} \\
  \;\;\;\;\;\;\;\; \;\;\;\;\;\;\; \;\;\;\;\;\;\; = \sum\limits_{i' = 1}^{{n_1}} {\sum\limits_{i = 1}^{{l_j}^\prime } {\left( {\begin{array}{*{20}{c}}
0\\
{{e_{i'}} \otimes y_i^{(j)}}
\end{array}} \right){{\left( {\begin{array}{*{20}{c}}
0\\
{{e_{i'}} \otimes y_i^{(j)}}
\end{array}} \right)}^T}} } \\
  \;\;\;\;\;\;\;\; \;\;\;\;\;\;\; \;\;\;\;\;\;\; = \sum\limits_{i' = 1}^{{n_1}} {\left( {\begin{array}{*{20}{c}}
0&0\\
0&{{e_{i'}}{{({e_{i'}})}^T} \otimes \sum\limits_{i = 1}^{{l_j}^\prime } {y_i^{(j)}{{(y_i^{(j)})}^T}} }
\end{array}} \right)} \\
 \;\;\;\;\;\;\;\; \;\;\;\;\;\;\; \;\;\;\;\;\;\; = \left( {\begin{array}{*{20}{c}}
0&0\\
0&{{I_{{n_1}}} \otimes {f_{{\eta _j}}}({G_2})}
\end{array}} \right),
\end{array}\]
where $j\in [q]\setminus\{1\}$. If $j=1$,
 $$F_{{\eta _1}}({G_1} \bigstar {G_2})={F_{k}}({G_1} \bigstar {G_2})=\left( {\begin{array}{*{20}{c}}
0&0\\
0&{{I_{{n_1}}} \otimes ({f_{{k}}}({G_2})-\frac{1}{n_{2}}J_{n_{2}})}
\end{array}} \right).
$$
Therefore, the spectral decomposition of $A(G_{1}\bigstar G_{2})$ is given by:
$$ A({G_1} \bigstar {G_2}) = \sum\limits_{i = 1}^p {\sum\limits_ \pm  {{\lambda _{i \pm }}{F_{{\lambda _{i \pm }}}}({G_1} \bigstar {G_2})} }  + \sum\limits_{j = 1}^q {{\eta _j}{F_{{\eta _j}}}} ({G_1} \bigstar {G_2}).$$

(ii) If there is some $i_{0}\in [p]$ such that $\lambda_{i_{0}}=0$, then $\lambda_{i_{0}+}=k$ and $\lambda_{i_{0}-}=0$. Similar to the proof of (i), we only need to compute $F_{k}(G_{1}\bigstar G_{2})$ and $F_{0}(G_{1}\bigstar G_{2})$. Let
$$ X_{{i_0} + }^j = \frac{1}{{{\sqrt{n_2}}}}\left( {\begin{array}{*{20}{c}}
0\\
{x_j^{({i_0})} \otimes {j_{{n_2}}}}
\end{array}} \right) \;
\text{and}\;
X_{{i_0} - }^j = \left( {\begin{array}{*{20}{c}}
{x_j^{({i_0})}}\\
0
\end{array}} \right)$$
for $j\in[l_{i_{0}}]$. Then
\[\small
A({G_1} \bigstar {G_2})X_{{i_0} + }^j = \frac{1}{{\sqrt{n_2}}}\left( {\begin{array}{*{20}{c}}
{A({G_1})}&{A({G_1}) \otimes j_{{n_2}}^T}\\
{A({G_1}) \otimes {j_{{n_2}}}}&{{I_{{n_1}}} \otimes A({G_2})}
\end{array}} \right)\left( {\begin{array}{*{20}{c}}
0\\
{x_j^{({i_0})} \otimes {j_{{n_2}}}}
\end{array}} \right) = \frac{k}{{\sqrt{n_2}}}\left( {\begin{array}{*{20}{c}}
0\\
{x_j^{({i_0})} \otimes {j_{{n_2}}}}
\end{array}} \right)
\]
and
\[
A({G_1} \bigstar {G_2})X_{{i_0} - }^j = \left( {\begin{array}{*{20}{c}}
{A({G_1})}&{A({G_1}) \otimes j_{{n_2}}^T}\\
{A({G_1}) \otimes {j_{{n_2}}}}&{{I_{{n_1}}} \otimes A({G_2})}
\end{array}} \right)\left( {\begin{array}{*{20}{c}}
{x_j^{({i_0})}}\\
0
\end{array}} \right) = 0\left( {\begin{array}{*{20}{c}}
{x_j^{({i_0})}}\\
0
\end{array}} \right)
\]
for $j\in [l_{i_{0}}]$. Thus, we have
\[
\begin{array}{l}
{F_k}({G_1} \bigstar {G_2}) = \sum\limits_{j = 1}^{{l_{_{{i_0}}}}} {X_{{i_0} + }^j{{(X_{{i_0} + }^j)}^T}}+\left( {\begin{array}{*{20}{c}}
0&0\\
0&{{I_{{n_1}}} \otimes {(f_{{k}}}({G_2})-\frac{1}{n_{2}}J_{n_{2}})}
\end{array}} \right) \\
\;\;\;\;\;\;\;\; \;\;\;\;\;\;\;\;\;\;\;\;= {\sum\limits_{j = 1}^{{l_{_{{i_0}}}}} {\frac{1}{{{n_2}}}\left( {\begin{array}{*{20}{c}}
0\\
{x_j^{({i_0})} \otimes {j_{{n_2}}}}
\end{array}} \right)\left( {\begin{array}{*{20}{c}}
0\\
{x_j^{({i_0})} \otimes {j_{{n_2}}}}
\end{array}} \right)} ^T}+\left( {\begin{array}{*{20}{c}}
0&0\\
0&{{I_{{n_1}}} \otimes {(f_{{k}}}({G_2})-\frac{1}{n_{2}}J_{n_{2}})}
\end{array}} \right)\\
\;\;\;\;\;\;\;\; \;\;\;\;\;\;\;\;\;\;\;\;= \frac{1}{{{n_2}}}\left( {\begin{array}{*{20}{c}}
0&0\\
0&{\sum\limits_{j = 1}^{{l_{_{{i_0}}}}} {x_j^{({i_0})}{{(x_j^{({i_0})})}^T} \otimes {J_{{n_2}}}} }
\end{array}} \right)+\left( {\begin{array}{*{20}{c}}
0&0\\
0&{{I_{{n_1}}} \otimes {(f_{{k}}}({G_2})-\frac{1}{n_{2}}J_{n_{2}})}
\end{array}} \right)\\
\;\;\;\;\;\;\;\; \;\;\;\;\;\;\;\;\;\;\;\;= \frac{1}{{{n_2}}}\left( {\begin{array}{*{20}{c}}
0&0\\
0&{{f_0}({G_1}) \otimes {J_{{n_2}}}}
\end{array}} \right)+\left( {\begin{array}{*{20}{c}}
0&0\\
0&{{I_{{n_1}}} \otimes {(f_{{k}}}({G_2})-\frac{1}{n_{2}}J_{n_{2}})}
\end{array}} \right)\\
\end{array}
\]
and
\[
{F_0}({G_1} \bigstar {G_2}) = \sum\limits_{j = 1}^{{l_{_{{i_0}}}}} {X_{{i_0} - }^j{{(X_{{i_0} - }^j)}^T}} = \left( {\begin{array}{*{20}{c}}
{\sum\limits_{j = 1}^{{l_{{i_{_0}}}}} {x_j^{({i_0})}{{(x_j^{({i_0})})}^T}} }&0\\
0&0
\end{array}} \right) = \left( {\begin{array}{*{20}{c}}
{{f_0}({G_1})}&0\\
0&0
\end{array}} \right).
\]
Therefore, the spectral decomposition of $A(G_{1}\bigstar G_{2})$ is given by:
\[
A({G_1} \bigstar {G_2}) = \sum\limits_{i \ne {i_0}}^{} {\sum\limits_ \pm  {{\lambda _{i \pm }}{F_{{\lambda _{i \pm }}}}({G_1} \bigstar {G_2})} }  + \sum\limits_{j = 2}^q {{\eta _j}{F_{{\eta _j}}}} ({G_1} \bigstar {G_2}) + k{F_k}({G_1} \bigstar {G_2}) + 0{F_0}({G_1} \bigstar {G_2}).
\]
 $\Box$\\
\\
\textbf{Proposition 3.3.} Let $G_{1}$ be any connected graph with $n_{1}$ vertices and $G_{2}$ be any $k$-regular graph with $n_{2}$ vertices. For $u$, $v\in V(G_{1})$,
\begin{enumerate}[(i)]
        \item if $\lambda_{j}\neq 0$ for any $j\in[p]$, then
        \begin{equation}
         {e_{(u,0)}}\exp ( - itA({G_1} \bigstar {G_2})){e_{(v,0)}} = \sum\limits_{j = 1}^p {{e^{ - it\frac{{{\lambda _j} + k}}{2}}}} {e_u}{f_{{\lambda _j}}}({G_1}){e_v}(\cos \frac{{{\Delta _{{\lambda _j}}}t}}{2} + i\frac{{k - {\lambda _j}}}{{{\Delta _{{\lambda _j}}}}}\sin \frac{{{\Delta _{{\lambda _j}}}t}}{2}),
        \end{equation}
         where $\Delta_{\lambda_{j}}=\sqrt{(\lambda_{j}-k)^{2}+4n_{2}\lambda_{j}^{2}}$ for $j\in [p]$.
        \item if there is some $j_{0}\in[p]$ such that $\lambda_{j_0}=0$, then
         \begin{equation}
         \begin{split}
         &{e_{(u,0)}}\exp ( - itA({G_1} \bigstar {G_2})){e_{(v,0)}}\\
         &\;\;\;\;\;\;\;\; \;\;\;\;\;\;= \sum\limits_{j \ne {j_0}}^{} {{e^{ - it\frac{{{\lambda _j} + k}}{2}}}} {e_u}{f_{{\lambda _j}}}({G_1}){e_v}(\cos \frac{{{\Delta _{{\lambda _j}}}t}}{2} + i\frac{{k - {\lambda _j}}}{{{\Delta _{{\lambda _j}}}}}\sin \frac{{{\Delta _{{\lambda _j}}}t}}{2}) + {e_u}{f_0}({G_1}){e_v},
         \end{split}
         \end{equation}
         where $\Delta_{\lambda_{j}}=\sqrt{(\lambda_{j}-k)^{2}+4n_{2}\lambda_{j}^{2}}$ for $j\in [p]$ and $j\neq j_{0}$.
     \end{enumerate}
\textit{Proof}: (i) Since $\lambda_{j}\neq0$ for any $j\in[p]$, according to (\ref{2}), then the transition matrix of $A(G_{1}\bigstar G_{2})$ is given by:
    \begin{equation}
\begin{array}{l}\label{4}
\exp ( - itA({G_1} \bigstar {G_2})) = \sum\limits_{j = 1}^p {\sum\limits_ \pm  {{e^{ - it{\lambda _{j \pm }}}}} } {F_{{\lambda _{j \pm }}}}({G_1} \bigstar {G_2}) + \sum\limits_{j' = 1}^q {{e^{ - it{\eta _{j'}}}}} {F_{{\eta _j}}}({G_1} \bigstar {G_2})\\
\end{array}.
\end{equation}
From (\ref{4}), (\ref{3}) and (\ref{5}), then the element of $\exp (- itA({G_1} \bigstar {G_2}))$ relevant to vertices $(u,0)$ and $(v,0)$ is given by:
\[
\begin{array}{l}
{e_{(u,0)}}\exp ( - itA({G_1} \bigstar {G_2})){e_{(v,0)}}\\
\\
 \;\;\;\;\;\;\;\; \;\;= {e_{(u,0)}}\sum\limits_{j = 1}^p {\sum\limits_ \pm  {{e^{ - it{\lambda _{j \pm }}}}} } {F_{{\lambda _{j \pm }}}}({G_1} \bigstar {G_2}){e_{(v,0)}} + {e_{(u,0)}}\sum\limits_{j' = 1}^p {{e^{ - it{\eta _{j'}}}}} {F_{{\eta _j}}}({G_1} \bigstar {G_2}){e_{(v,0)}}\\
\\
 \;\;\;\;\;\;\;\; \;\;= \sum\limits_{j = 1}^p {{e^{ - it\frac{{{\lambda _j} + k}}{2}}}} \sum\limits_ \pm  {{e^{ \mp it\frac{\Delta _{\lambda _j}}{2}}}} \frac{{{\lambda _j}^2}}{{{{({\lambda _{j \pm }} - k)}^2} + {n_2}{\lambda _j}^2}}\frac{{{{({\lambda _{j \pm }} - k)}^2}}}{{{\lambda _j}^2}}{e_u}{f_{{\lambda _j}}}({G_1}){e_v}\\
\\
 \;\;\;\;\;\;\;\; \;\;= \sum\limits_{j = 1}^p {{e^{ - it\frac{{{\lambda _j} + k}}{2}}}} {e_u}{f_{{\lambda _j}}}({G_1}){e_v}\sum\limits_ \pm  {{e^{ \mp it\frac{{{\Delta _{{\lambda _j}}}}}{2}}}} \frac{{{{({\lambda _{j \pm }} - k)}^2}}}{{{{({\lambda _{j \pm }} - k)}^2} + {n_2}{\lambda _j}^2}}\\
\\
 \;\;\;\;\;\;\;\; \;\;= \sum\limits_{j = 1}^p {{e^{ - it\frac{{{\lambda _j} + k}}{2}}}} {e_u}{f_{{\lambda _j}}}({G_1}){e_v}(\frac{{{{({\lambda _{j + }} - k)}^2}}}{{{{({\lambda _{j + }} - k)}^2} + {n_2}{\lambda _j}^2}} + \frac{{{{({\lambda _{j - }} - k)}^2}}}{{{{({\lambda _{j - }} - k)}^2} + {n_2}{\lambda _j}^2}})\cos \frac{{{\Delta _{{\lambda _j}}}t}}{2}\\
 \\
\;\;\;\;\;\;\;\; \;\;\;\;\; + (\frac{{{{({\lambda _{j - }} - k)}^2}}}{{{{({\lambda _{j - }} - k)}^2} + {n_2}{\lambda _j}^2}} - \frac{{{{({\lambda _{j + }} - k)}^2}}}{{{{({\lambda _{j + }} - k)}^2} + {n_2}{\lambda _j}^2}})i\sin \frac{{{\Delta _{{\lambda _j}}}t}}{2}.
\end{array}
\]
In the light of (\ref{6}) and (\ref{7}), we obtain the following three equalities,
\begin{equation}\label{8}
(\lambda_{j+}-k)(\lambda_{j-}-k)=-n_{2}\lambda_{j}^{2},
\end{equation}
\begin{equation}\label{9}
(\lambda_{j+}-k)^{2}+(\lambda_{j-}-k)^{2}=(\lambda_{j}-k)^{2}+2n_{2}\lambda_{j}^{2},
\end{equation}
and
\begin{equation}\label{9*}
(\lambda_{j-}-k)^{2}-(\lambda_{j+}-k)^{2}=\Delta_{\lambda_{j}}(k-\lambda_{j}).
\end{equation}
Thus,
$$\frac{{{{({\lambda _{j + }} - k)}^2}}}{{{{({\lambda _{j + }} - k)}^2} + {n_2}{\lambda _j}^2}} + \frac{{{{({\lambda _{j - }} - k)}^2}}}{{{{({\lambda _{j - }} - k)}^2} + {n_2}{\lambda _j}^2}}=1$$
and
$$\frac{{{{({\lambda _{j - }} - k)}^2}}}{{{{({\lambda _{j - }} - k)}^2} + {n_2}{\lambda _j}^2}} - \frac{{{{({\lambda _{j + }} - k)}^2}}}{{{{({\lambda _{j + }} - k)}^2} + {n_2}{\lambda _j}^2}}=\frac{k-\lambda_{j}}{\Delta_{\lambda_{j}}}.$$
Therefore,
$${e_{(u,0)}}\exp ( - itA({G_1} \bigstar {G_2})){e_{(v,0)}} = \sum\limits_{j = 1}^p {{e^{ - it\frac{{{\lambda _j} + k}}{2}}}} {e_u}{f_{{\lambda _j}}}({G_1}){e_v}(\cos \frac{{{\Delta _{{\lambda _j}}}t}}{2} + i\frac{{k - {\lambda _j}}}{{{\Delta _{{\lambda _j}}}}}\sin \frac{{{\Delta _{{\lambda _j}}}t}}{2}).$$

(ii) If there is some $j_{0}\in[p]$ such that $\lambda_{j_{0}}=0$, then the proof is totally similar to that of (i), here we omit it. $\Box$

\section{Perfect State Transfer}

In this section, we mainly focus on PST of neighborhood corona of two graphs. First we recall the following lemma.\\
\\
\textbf{Lemma 4.1.} (Li, Liu and Zhang \cite{edge corona}) Let $G$ be a $k$-regular connected graph with $n$ vertices, for any vertex $v$ of $G$.
\begin{enumerate}[(i)]
\item If $G$ is not a complete graph, then $|\text{supp}_{G}(v)|\geq 3$;
\item If $G$ is a complete graph, then $|\text{supp}_{G}(v)|=2$.
\end{enumerate}
\textit{ Proof}: Suppose that $\{\lambda_{i}|i\in[p]\}$ is the set of all distinct eigenvalues of $A(G)$ and $\lambda_{1}> \lambda_{2}>\cdots>\lambda_{p}$. Denote the eigenprojector of eigenvalue $\lambda_{i}$ by $f_{\lambda_{i}}(G)$. Since $G$ is $k$-regular, then $L(G)=D(G)-A(G)=kI_{n}-A(G)$. Let $\{\theta_{i}|i\in[p]\}$ be the set of all distinct eigenvalues of $L(G)$ and $\theta_{1}< \theta_{2}<\cdots<\theta_{p}$, then $\theta_{i}=k-\lambda_{i}$. It is clear that $f_{\lambda_{i}}(G)=f_{\theta_{i}}(G)$ for the regular graph $G$. Since, for any $\theta_{i}\in \text{supp}_{L(G)}(v)$, $f_{\theta_{i}}(G)e_{v}\neq 0$, then $\lambda_{i}\in\text{supp}_{G}(v)$ for $f_{\lambda_{i}}(G)e_{v}=f_{\theta_{i}}(G)e_{v}\neq0 $. Conversely, since, for any $\lambda_{i}\in\text{supp}_{G}(v)$, $f_{\lambda_{i}}(G)e_{v}\neq 0$, then $\theta_{i}\in\text{supp}_{L(G)}(v)$ for $f_{\theta_{i}}(G)e_{v}=f_{\lambda_{i}}(G)e_{v}\neq0$. Therefore, we obtain the desired results from Lemma 4.2 in \cite{edge corona}.\\
\\
\textbf{Lemma 4.2.} If $0$ is not an eigenvalue of $A(G_{1})$ and $G_{2}$ is a $k$-regular graph. Then $(v, 0)$ is periodic whenever $(v, w)$ is periodic in $G_{1}\bigstar G_{2}$ for $v\in V(G_{1})$, $w\in V(G_{2})$.\\
\\
\textit{ Proof}: According to (\ref{3}) and (\ref{5}), it is obvious that the element of $\text{supp}_{G_{1}\bigstar {G_{2}}}(v, 0)$ contains in $\text{supp}_{G_{1}\bigstar {G_{2}}}(v, w)$. If $(v, w)$ is periodic in $G_{1}\bigstar G_{2}$, then $(v, 0)$ is periodic in $G_{1}\bigstar G_{2}$  according to Lemma 2.3.\\
\\
\textbf{Lemma 4.3.} Let $G_{1}$ be a connected graph and $G_{2}$ be a $k$-regular graph. Suppose that there is some $i\in [p]$ such that the eigenvalue $\lambda_{i}=0$ for $A(G_{1})$. For any $v\in V(G_{1})$ and $w\in V(G_{2})$,
\begin{enumerate}[(i)]
\item if $0\notin\text{supp}_{G_{1}}(v)$, then $(v, 0)$ is periodic in $G_{1}\bigstar G_{2}$ whenever $(v, w)$ is periodic in $G_{1}\bigstar G_{2}$;
\item if $0\in\text{supp}_{G_{1}}(v)$, then $\text{supp}_{G_{1}\bigstar {G_{2}}}(v, 0)\setminus\{0\} \subseteq \text{supp}_{G_{1}\bigstar {G_{2}}}(v, w)$.
\end{enumerate}
\textit{Proof}: (i) Assume that $0\notin\text{supp}_{G_{1}}(v)$. According to the (ii) of Proposition 3.2, one has $\text{supp}_{G_{1}\bigstar {G_{2}}}(v, 0) \subseteq \text{supp}_{G_{1}\bigstar {G_{2}}}(v, w)$. Now, we obtain the desired results from Lemma 2.3.

(ii) If $0\in\text{supp}_{G_{1}}(v)$, then $f_{0}(G_{1})e_{v}\neq 0$. According to (\ref{9}), we easily obtain that $F_{0}(G_{1}\bigstar {G_{2}})e_{(v,0)}\neq0$ and  $F_{0}(G_{1}\bigstar {G_{2}})e_{(v,w)}=0$. Hence, $0\in \text{supp}_{G_{1}\bigstar {G_{2}}}(v, 0)$, but $0\notin \text{supp}_{G_{1}\bigstar {G_{2}}}(v, w)$.\\
\\
\textbf{Theorem 4.4.} Assume that $G_{1}$ is an $r$-regular connected integral graph with $n_{1}$ vertices. Let $G_{2}$ be a $k$-regular graph with $n_{2}$ vertices. If $\sqrt{(r-k)^{2}+4n_{2}r^{2}}$ is not an integer, then $(v, w)$ is not periodic for any $v\in V(G_{1})$, $w\in V(G_{2})\cup\{0\}$ in $G_{1}\bigstar G_{2}$. Moreover, there is no PST in $G_{1}\bigstar G_{2}$.\\
\\
\textit{Proof}: Since $G_{1}$ is $r$-regular, then $A(G_{1})j_{n_{1}}=rj_{n_{1}}$. According to the definition of eigenprojector, the eigenprojector of eigenvalue $r$ in $G_{1}$ is given by $f_{r}(G_{1})=\frac{1}{n_{1}}J_{n_{1}}$ as $G_{1}$ is a connected graph. Therefore, $f_{r}(G_{1})e_{v}\neq0$ for any $v\in V(G_{1})$. In other words, $r\in\text{supp}_{G_{1}}(v)$ for any $v\in V(G_{1})$.

Suppose that $G_{1}\bigstar G_{2}$ admits PST between vertex $(v, 0)$ and another vertex.  Then $(v, 0)$ is periodic in $G_{1}\bigstar G_{2}$. According to Lemma 2.3 and Theorem 2.1, non-zero elements of $\text{supp}_{G_{1}\bigstar G_{2}}(v, 0)$ are all integers or the form of $\frac{a+b\sqrt{\Delta}}{2}$ for integer $a$, square-free integer $\Delta$ and some integer $b$.

Case $1$: Non-zero elements of $\text{supp}_{G_{1}\bigstar G_{2}}(v, 0)$ are all integers.
Since $r\in\text{supp}_{G_{1}}(v)$, then $F_{r\pm}(G_{1}\bigstar G_{2})e_{(v, 0)}\neq0$ for $f_{r}(G_{1})e_{v}\neq0$. Hence $\lambda_{r\pm}\in\text{supp}_{G_{1}\bigstar G_{2}}(v, 0)$. Since all the elements of $\text{supp}_{G_{1}\bigstar G_{2}}(v, 0)$ are integers, then $\lambda_{r+}$, $\lambda_{r-}$, and $\lambda_{r+}-\lambda_{r-}$ are integers.
According to  $\lambda_{r+}=\frac{r+k+\sqrt{(r-k)^{2}+4n_{2}r^{2}}}{2}$ and $\lambda_{r-}=\frac{r+k-\sqrt{(r-k)^{2}+4n_{2}r^{2}}}{2}$, one has
\begin{equation}\label{10}
\lambda_{r+}-\lambda_{r-}=\sqrt{(r-k)^{2}+4n_{2}r^{2}}.
\end{equation}
The left of (\ref{10}) is an integer but the right of (\ref{10}) is not an integer. This contradicts that non-zero elements of $\text{supp}_{G_{1}\bigstar G_{2}}(v, 0)$ are integers.

Case $2$: All the elements of $\text{supp}_{G_{1}\bigstar G_{2}}(v, 0)$ are the form of $\frac{a+b\sqrt{\Delta}}{2}$ for integer $a$, square-free integer $\Delta$ and some integer $b$. According to Lemma 4.1, there is $\lambda\in\text{supp}_{G_{1}}(v)$ and  $\lambda\neq r$, such that $\lambda_{\pm}\in\text{supp}_{G_{1}\bigstar G_{2}}(v, 0)$ for $f_{\lambda}(G_{1})e_{v}\neq0$. Notice that $r_{\pm}\in\text{supp}_{G_{1}\bigstar G_{2}}(v, 0)$ and $\sqrt{(r-k)^{2}+4n_{2}r^{2}}$ is not an integer. Let $\lambda_{\pm}=\frac{a+b_{\pm}\sqrt{\Delta}}{2}$ for integer $a$, square-free integer $\Delta>1$ and some integer $b_{\pm}$. In the light of (\ref{8}), one has
\begin{equation}\label{11}
(\lambda_{+}-k)(\lambda_{-}-k)=-n_{2}\lambda^{2}.
\end{equation}
According to the form of $\lambda_{\pm}$, we obtain the following equality:
\begin{equation}\label{12}
\frac{1}{4}[(a-2k)^{2}+b_{+}b_{-}\Delta]+\frac{1}{4}(a-2k)(b_{+}+b_{-})\sqrt{\Delta}=-n_{2}\lambda^{2}.
\end{equation}
Since $\sqrt{\Delta}$ is irrational, (\ref{12}) holds if and only if either $a-2k=0$, or $b_{+}+b_{-}=0$.
If $b_{+}+b_{-}=0$, then $a=\lambda_{+}+\lambda_{-}=\lambda+k$, which implies that $|\text{supp}_{G_{1}}(v)|=1$. At that time, it contradicts with Lemma 4.1. If $a-2k=0$, then $\lambda_{\pm}=k+\frac{b_{\pm}\sqrt{\Delta}}{2}$. Therefore,
\begin{equation}\label{13}
(\lambda_{+}-k)+(\lambda_{-}-k)=\lambda-k.
\end{equation}
The left of (\ref{13}) is a rational multiple of $\sqrt{\Delta}$, but the right of (\ref{13}) is an integer. Thus, (\ref{13}) can not hold, which means that $\lambda_{\pm}$ can not be the form of quadratic integer. Therefore, not all the elements of $\text{supp}_{G_{1}\bigstar G_{2}}(v, 0)$ are the form of quadratic integer.

According to Case 1 and Case 2, $(v, 0)$ is not periodic in $G_{1}\bigstar G_{2}$ for any $v\in V(G_{1})$ by Lemma 2.3. Furthermore, $(v, w)$ is not periodic in $G_{1}\bigstar G_{2}$ for any $v\in V(G_{1})$ and any $w\in V(G_{2})\cup\{0\}$ by Lemma 4.2. Finally, according to Lemma 2.2, there is no PST in $G_{1}\bigstar G_{2}$.
$\Box$\\

Remark that Theorem 4.4 shows that there is no PST in $G_{1}\bigstar G_{2}$ when $\sqrt{(r-k)^{2}+4n_{2}r^{2}}$ is not an integer for two regular graphs $G_{1}$ and $G_{2}$. Next, we consider the case of non-regular graphs.\\
\\
\textbf{Theorem 4.5.} Assume that $G_{1}$ is any connected graph. Let $G_{2}$ be any connected $k$-regular graph with $m_{2}$ vertices. If there exists $\lambda_{i}\in \text{supp}_{G_{1}}(v)$ such that $\lambda_{i}\in \mathbb{Z}\sqrt{\Delta'}$ for some $i\in [p]$ and some square-free integer $\Delta'>1$, where $\mathbb{Z}$ is the set of all integers. Then $(v, w)$ is not periodic in $G_{1}\bigstar G_{2}$ for $v\in V(G_{1})$, and any $w\in V(G_{2})\cup\{0\}$.\\
\\
\textit{Proof}: Suppose that there is PST in $G_{1}\bigstar G_{2}$ between vertex $(v, 0)$ and another vertex, then $(v, 0)$ is periodic in $G_{1}\bigstar G_{2}$. Thus, according to Lemma 2.3 and Theorem 2.1, non-zero elements of $\text{supp}_{G_{1}\bigstar G_{2}}(v, 0)$ are either all integers, or all quadratic integers.

Case $1$: Non-zero elements of $\text{supp}_{G_{1}\bigstar G_{2}}(v, 0)$ are all integers.

Since $\lambda_{i}\in \text{supp}_{G_{1}}(v)$, then $f_{\lambda_{i}}(G_{1})e_{v}\neq 0$. According to (\ref{3}), then $\lambda_{i\pm}\in\text{supp}_{G_{1}\bigstar G_{2}}(v, 0)$ for $F_{\lambda_{i\pm}}(G_{1}\bigstar G_{2})e_{(v, 0)}\neq 0$. Hence $\lambda_{i+}$ and $\lambda_{i-}$ are both integers. Furthermore, $\lambda_{i+}+ \lambda_{i-}$ is also an integer. According to (\ref{6}) and (\ref{7}),
\begin{equation}\label{14}
\lambda_{i+}+ \lambda_{i-}=\lambda_{i} +k.
\end{equation}
Since $\lambda_{i}\in \mathbb{Z}\sqrt{\Delta'}$ is an irrational number and the left of equality (\ref{14}) is an integer. However, the right of (\ref{14}) is an irrational number, it is impossible. Hence not all $\text{supp}_{G_{1}\bigstar G_{2}}(v, 0)$ are integers.

Case $2$: There is some integer $a$ and square-free integer $\Delta$ such that all the elements of $\text{supp}_{G_{1}\bigstar G_{2}}(v, 0)$ are the form of $\frac{a+b_{\lambda_{i\pm}}\sqrt{\Delta}}{2}$ for some integer $b_{\lambda_{i\pm}}$ relevant to $\lambda_{i\pm}$.

According to Case 1, $\lambda_{i\pm}\in\text{supp}_{G_{1}\bigstar G_{2}}(v, 0)$. Since $(\lambda_{i}-k)^{2}+4n_{2}\lambda_{i}^{2}$ is an irrational number according to $\lambda_{i}\in Z\sqrt{\Delta'}$, then there is an integer $a$ and square-free integer $\Delta>1$  such that $\lambda_{i\pm}=\frac{a+b_{\lambda_{i\pm}}\sqrt{\Delta}}{2}$ for some integer $b_{\lambda_{i\pm}}$ corresponding to $\lambda_{i\pm}$.
In the light of (\ref{8}), then $(\lambda_{i+}-k)(\lambda_{i-}-k)=-n_{2}\lambda_{i}^{2}$. According to the form of $\lambda_{i\pm}$, one has, by a simple calculation,
\begin{equation}\label{15}
\frac{1}{4}[(a-2k)^{2}+b_{\lambda_{i+}}b_{\lambda_{i-}}\Delta]+\frac{1}{4}(a-2k)(b_{\lambda_{i+}}+b_{\lambda_{i-}})\sqrt{\Delta}=-n_{2}\lambda_{i}^{2}.
\end{equation}
Since $\lambda_{i}\in \mathbb{Z}\sqrt{\Delta'}$, then $-n_{2}\lambda_{i}^{2}$ is an integer. Since $\Delta>1$, then $\sqrt{\Delta}$ is an irrational number. In the light of (\ref{15}), we obtain that either $a-2k=0$, or $b_{\lambda_{i+}}+b_{\lambda_{i-}}=0$. If $b_{\lambda_{i+}}+b_{\lambda_{i-}}=0$, then $a=\lambda_{i+}+\lambda_{i-}=\lambda_{i}+k$. It contradicts that $a$ is an integer beacause $\lambda_{i}\in \mathbb{Z}\sqrt{\Delta'}$ is an irrational number. If $a-2k=0$, then $a=2k$. Hence $\lambda_{i+}=k+\frac{b_{\lambda_{i+}}\sqrt{\Delta}}{2}$ and $\lambda_{i-}=k+\frac{b_{\lambda_{i-}}\sqrt{\Delta}}{2}$. According to the equalities above, then
\begin{equation}\label{16}
(\lambda_{i+}-k)+(\lambda_{i-}-k)=\lambda_{i}-k
\end{equation}
Taking square both sides of (\ref{16}), we easily obtain that the left is a rational number, but the right is irrational, a contradiction. Therefore, neither $a-2k=0$ nor $b_{\lambda_{i+}}+b_{\lambda_{i-}}=0$. Furthermore, non-zero elements of  $\text{supp}_{G_{1}\bigstar G_{2}}(v, 0)$ are all quadratic integers.

Now, it follows from Case $1$, Case $2$ and Lemma 2.3 that $(v, 0)$ is not periodic in $G_{1}\bigstar G_{2}$. Therefore, $(v, w)$ is not periodic in $G_{1}\bigstar G_{2}$. $\Box$\\
\\
\textbf{Theorem 4.6.} The neighborhood corona $P_{3}\bigstar G$ has no PST for a connected $k$-regular graph $G$.\\
\\
\textit{Proof}: It is known that $\text{Sp}(P_{3})=\{\sqrt{2},0,-\sqrt{2}\}$, By a simple calculation, then the eigenvalue support of middle vertex is the set $\{\sqrt{2},-\sqrt{2}\}$ and the eigenvalue support the endpoints of $P_{3}$ are the set $\{\sqrt{2},0,-\sqrt{2}\}$. Based on Theorem 4.5, $(v,0)$ is not a periodic vertex in $P_{3}\bigstar G$ for any $v\in V(P_{3})$, then $(v,w)$ is not a periodic vertex in $P_{3}\bigstar G$ for any $v\in V(P_{3})$, $w\in V(G)\cup\{0\}$. According to Lemma 2.2, $P_{3}\bigstar G$ has no PST.   $\Box$

\section{Pretty Good State Transfer}

Since the neighborhood corona graphs having PST are rare. Then it is meaningful for us to search some neighborhood corona graphs admitting PGST. In this section, we mainly study PGST in neighborhood corona of two graphs.\\
\\
\textbf{Theorem 5.1.} Assume that $G_{1}$ has PST at time $\frac{\pi}{g}$ between vertices $x$ and $y$, where $g$ is defined in Theorem 2.1. If $\sqrt{4n_{2}+1}$ is not an integer, then there is PGST in $G_{1}\bigstar \overline{K_{n_{2}}}$ between vertices $(x, 0)$ and $(y, 0)$.\\
\\
\textit{Proof}: Since $G_{1}$ has PST at time $\frac{\pi}{g}$, then all eigenvalues $\lambda_{j}$ are integers by Theorem 2.1 for $j\in[p]$.
The following proof is divided into two cases. \\

Case $1$: $0\notin \text{supp}_{G_{1}}(y)$. \\

According to (i) of Proposition 3.3, we have
\begin{equation}
{e_{(x,0)}}\exp ( - itA({G_1} \bigstar \overline{K_{n_{2}}})){e_{(y,0)}} = \sum\limits_{j = 1}^p {{e^{ - it\frac{{{\lambda _j}}}{2}}}} {e_x}{f_{{\lambda _j}}}({G_1}){e_y}(\cos \frac{{{\Delta _{{\lambda _j}}}t}}{2} + i\frac{-{\lambda _j}}{{{\Delta _{{\lambda _j}}}}}\sin \frac{{{\Delta _{{\lambda _j}}}t}}{2}),
\end{equation}
where $\Delta _{{\lambda _j}}=\sqrt{(4n_{2}+1)\lambda_{j}^{2}}$ for $j\in[p]$.
Since $f_{{\lambda _j}}({G_1}){e_y}=0$ for $\lambda _j \notin \text{supp}_{G_{1}}(y)$, then
\begin{equation}\small
 {e_{(x,0)}}\exp ( - itA({G_1} \bigstar \overline{K_{n_{2}}})){e_{(y,0)}}=\sum\limits_{\lambda _j\in \text{supp}_{G_{1}}(y)} {{e^{ - it\frac{{{\lambda _j}}}{2}}}} {e_x}{f_{{\lambda _j}}}({G_1}){e_y}(\cos \frac{{{\Delta _{{\lambda _j}}}t}}{2} + i\frac{{- {\lambda _j}}}{{{\Delta _{{\lambda _j}}}}}\sin \frac{{{\Delta _{{\lambda _j}}}t}}{2}).
\end{equation}
In order to prove that ${G_1} \bigstar \overline{K_{n_{2}}}$ has PGST between $(x, 0)$ and $(y, 0)$, we need to find a time $t_{0}$ such that
 $$\left| \sum\limits_{\lambda _j\in \text{supp}_{G_{1}}(y)} {{e^{ - it_{0}\frac{{{\lambda _j}}}{2}}}} {e_x}{f_{{\lambda _j}}}({G_1}){e_y}(\cos \frac{{{\Delta _{{\lambda _j}}}t_{0}}}{2} + i\frac{{- {\lambda _j}}}{{{\Delta _{{\lambda _j}}}}}\sin \frac{{{\Delta _{{\lambda _j}}}t_{0}}}{2})\right|\approx 1.$$
Since $\sqrt{4n_{2}+1}$ is not an integer, then $\Delta _{{\lambda _j}}=\sqrt{(4n_{2}+1)\lambda_{j}^{2}}$ is not an integer for integer $\lambda_{j}\in \text{supp}_{G_{1}}(y)$.
Let $\Delta _{{\lambda _j}}=a_{j}\sqrt{b_{j}}$ for each $\lambda_{j}\in \text{supp}_{G_{1}}(y)$, where $a_{j}$, $ b_{j} \in \mathbb{Z}^{+}$ and $b_{j}$ is the square-free part of $\Delta _{{\lambda _j}}^{2}$.
Then the disjoint union $\{1\}\cup \{\sqrt{b_{j}}:\lambda _j\in \text{supp}_{G_{1}}(y)\}$ is linearly independent over $\mathbb{Q}$ by Lemma 2.5. According to Theorem 2.4, there are integers $\alpha$ and $c_{j}$ for each $\lambda_{j}\in \text{supp}_{G_{1}}(y)$ such that
\begin{equation}\label{17}
\alpha \sqrt{b_{j}}-c_{j}\approx-\frac{1}{2g}\sqrt{b_{j}}.
\end{equation}
If $b_{j}=b_{j'}$ for two different eigenvalues $\lambda_{j}$, $\lambda_{j'}\in \text{supp}_{G_{1}}(y)$, then $c_{j}=c_{j'}$. Multiplying $4a_{j}$ to two sides of (\ref{17}), we have
$\Delta _{{\lambda _j}}\approx\frac{4a_{j}c_{j}}{4\alpha+\frac{2}{g}}$.
Let $t_{0}=(4\alpha+\frac{2}{g})\pi $. Then $\cos \frac{\Delta _{{\lambda _j}}t_{0}}{2}\approx \cos2a_{j}c_{j}\pi=1$.
Hence,
$$\begin{array}{l}
 \left|{e_{(x,0)}}\exp ( - it_{0}A({G_1} \bigstar \overline{K_{n_{2}}})){e_{(y,0)}}\right|
 =\left|\sum\limits_{\lambda _j\in \text{supp}_{G_{1}}(y)} {{e^{ - it_{0}\frac{{{\lambda _j}}}{2}}}} {e_x}{f_{{\lambda _j}}}({G_1}){e_y}(\cos \frac{{{\Delta _{{\lambda _j}}}t_{0}}}{2} + i\frac{{- {\lambda _j}}}{{{\Delta _{{\lambda _j}}}}}\sin \frac{{{\Delta _{{\lambda _j}}}t_{0}}}{2})\right|\\
 \\
 \;\;\;\;\;\;\;\;\;\;\;\;\;\;\;\;\;\;\;\;\;\;\;\;\;\;\;\;\;\;\;\;\;\;\;\;\;\;\;\;\;\;\;\;\;\;\;\;\;\;\;\;\;\;\approx\left |\sum\limits_{\lambda _j\in \text{supp}_{G_{1}}(y)} {{e^{ - it_{0}\frac{{{\lambda _j}}}{2}}}} {e_x}{f_{{\lambda _j}}}({G_1}){e_y}\right|\\
 \\
 \;\;\;\;\;\;\;\;\;\;\;\;\;\;\;\;\;\;\;\;\;\;\;\;\;\;\;\;\;\;\;\;\;\;\;\;\;\;\;\;\;\;\;\;\;\;\;\;\;\;\;\;\;\;=\left|\sum\limits_{\lambda _j\in \text{supp}_{G_{1}}(y)} {{e^{ - i\frac{\pi}{g}\lambda_{j}}}} {e_x}{f_{{\lambda _j}}}({G_1}){e_y}\right|\\
 \\
 \;\;\;\;\;\;\;\;\;\;\;\;\;\;\;\;\;\;\;\;\;\;\;\;\;\;\;\;\;\;\;\;\;\;\;\;\;\;\;\;\;\;\;\;\;\;\;\;\;\;\;\;\;\;=1,
\end{array}$$
where the last equality holds by Theorem 2.1. \\

Case $2$: $0\in \text{supp}_{G_{1}}(y)$.\\

According to (ii) of Proposition 3.3, if there exists $0=\lambda _{j_{0}}\in\text{supp}_{G_{1}}(y)$, then
\begin{equation}\small
{e_{(x,0)}}\exp ( - itA({G_1} \bigstar \overline{K_{n_{2}}})){e_{(y,0)}} = \sum\limits_{j\neq j_{0}} {{e^{ - it\frac{{{\lambda _j}}}{2}}}} {e_x}{f_{{\lambda _j}}}({G_1}){e_y}(\cos \frac{{{\Delta _{{\lambda _j}}}t}}{2} + i\frac{-{\lambda _j}}{{{\Delta _{{\lambda _j}}}}}\sin \frac{{{\Delta _{{\lambda _j}}}t}}{2})+{e_x}{f_{0}}({G_1}){e_y},
\end{equation}
where $\Delta _{{\lambda _j}}=\sqrt{(4n_{2}+1)\lambda_{j}^{2}}$ for $j\neq j_{0}$.
Similar to the discussion in Case $1$, we can obtain desired result based on the equalities $1=e^{-i\frac{\pi}{g}0}$ and $e^{ - it_{0}\frac{{{\lambda _j}}}{2}}=e^{-i\frac{\pi}{g}\lambda_{j}}$.

Therefore, there exists PGST on ${G_1} \bigstar \overline{K_{n_{2}}}$ between $(x, 0)$ and $(y, 0)$. $\Box$\\

Theorem 5.1 gives a sufficient condition such that $G_{1}\bigstar \overline{K_{n_{2}}}$ admits PGST when $G_{1}$ has PST at $\frac{\pi}{g}$. It is well known that $C_{4}$ has PST between its two antipodal vertices at $\frac{\pi}{2}$. Therefore,  we immediately obtain the following corollary.\\
\\
\textbf{Corollary 5.2.} If $\sqrt{1+4n_{2}}$ is not an integer, then ${C_{4}} \bigstar \overline{K_{n_{2}}}$ admits PGST.\\
\\
\textit{Proof}: Recall that $\text{Sp}(C_{4})=\{2, 0^{2}, -2\}$ and $C_{4}$ has PST at time $\frac{\pi}{2}$ between two antipodal vertices (see \cite{Coutinho2014}). Without loss of generality, suppose that $C_{4}$ admits PST between vertices $v_{1}$ and $v_{3}$. According to the definition of eigenvalue support, then $\text{supp}_{C_{4}}{(v_{1})}=\{2, 0, -2\}=\text{supp}_{C_{4}}{(v_{3})}$.
Since $\sqrt{1+4n_{2}}$ is not an integer and $e^{-i\frac{\pi}{2}0}=1$.
Therefore, ${C_{4}} \bigstar \overline{K_{n_{2}}}$ admits PGST at vertices $(v_{1},0)$ and $(v_{3},0)$ by Theorem 5.1.  $\Box$\\

Theorem 5.2 implies that ${C_{4}} \bigstar \overline{K_{n_{2}}}$ admits PGST when $\sqrt{1+4n_{2}}$ is not an integer and $G_{2}=\overline{K_{n_{2}}}$ is not a connected graph. Next, we consider PGST on ${C_4} \bigstar G_{2}$ whenever $G_{2}$ is a $k$-regular connected graph.\\
\\
\textbf{Theorem 5.3.} Assume that $G_{2}$ is a $k$-regular connected graph with vertices $n_2$. If neither of $\sqrt{(2+k)^{2}+16n_{2}}$ and $\sqrt{(2-k)^{2}+16n_{2}}$ is an integer and $k=0(\text{mod} 4)$, then ${C_4} \bigstar G_{2}$ admits PGST.\\
\\
\textit{Proof}: Suppose that $V(C_{4})=\{v_{1}, v_{2}, v_{3}, v_{4}\}$. According to the proof of Theorem 5.2, $\text{Sp}(C_{4})=\{2, 0^{2}, -2\}$ and $\text{supp}_{C_{4}}(v_{1})=\text{supp}_{C_{4}}(v_{3})=\{2,0,-2\}$. Let $S=\text{supp}_{C_{4}}(v_{3})$. Since ${f_{{\lambda _j}}}({C_4}){e_{v_{3}}}=0$ for any $\lambda_{j}\notin S$, according to Proposition 3.3,
 \begin{equation}
 \begin{split}
 &{e_{(v_{1},0)}}\exp ( - itA({C_{4}} \bigstar {G_{2}})){e_{(v_{3},0)}}\\
 &\;\;\;\;\;\;\;= \sum\limits_{\lambda_{j}\neq0} {{e^{ - it\frac{{{\lambda _j} + k}}{2}}}} {e_{v_{1}}}{f_{{\lambda _j}}}({C_4}){e_{v_{3}}}(\cos \frac{{{\Delta _{{\lambda _j}}}t}}{2} + i\frac{{k - {\lambda _j}}}{{{\Delta _{{\lambda _j}}}}}\sin \frac{{{\Delta _{{\lambda _j}}}t}}{2}) + {e_{v_{1}}}{f_0}({C_4}){e_{v_{3}}},
 \end{split}
 \end{equation}
 where $\Delta_{\lambda_{j}}=\sqrt{(\lambda_{j}-k)^{2}+4n_{2}\lambda_{j}^{2}}$ for $j\in\{1, 3\}$. Hence, in order to prove that ${C_4} \bigstar G_{2}$ occurs PGST between vertices $(v_{1},0)$ and $(v_{3},0)$, we only need to find some $t_{0}$ such that
\begin{equation}
\begin{split}
 &\left|{e_{(v_{1},0)}}\exp ( - it_{0}A({C_{4}} \bigstar {G_{2}})){e_{(v_{3},0)}}\right|\\
 &\;\;\;\;\;= \left|\sum\limits_{\lambda_{j}\neq0} {{e^{ - it_{0}\frac{{{\lambda _j} + k}}{2}}}} {e_{v_{1}}}{f_{{\lambda _j}}}({C_4}){e_{v_{3}}}(\cos \frac{{{\Delta _{{\lambda _j}}}t_{0}}}{2} + i\frac{{k - {\lambda _j}}}{{{\Delta _{{\lambda _j}}}}}\sin \frac{{{\Delta _{{\lambda _j}}}t_{0}}}{2}) + {e_{v_{1}}}{f_0}({C_4}){e_{v_{3}}}\right|\\
 &\;\;\;\;\;\approx 1.
\end{split}
\end{equation}
Since $\lambda_{1}=2$ and $\lambda_{3}=-2$, obviously, $\Delta_{\lambda_{3}}=\sqrt{(2+k)^{2}+16n_{2}}$ and $\Delta_{\lambda_{1}}=\sqrt{(2-k)^{2}+16n_{2}}$.
Based on the given condition, neither of $\sqrt{(2+k)^{2}+16n_{2}}$ and $\sqrt{(2-k)^{2}+16n_{2}}$ is an integer, then $\Delta_{\lambda_{j}}=a_{j}\sqrt{b_{j}}$ for $j\in\{1, 3\}$, where $a_{j}, b_{j}\in \mathbb{Z}^{+}$ and $b_{j}$ is square-free part of $\Delta_{\lambda_{j}}^{2}$.
It is not difficult to see that $\{1\}\cup\{\sqrt{b_{j}}: j\in\{1, 3\}\}$ is linearly independent over $\mathbb{Q}$.
According to Theorem 2.4, there are integers $\alpha$ and $d_{j}$ such that
\begin{equation}\label{18}
\alpha \sqrt{b_{j}}-d_{j}\approx-\frac{1}{4}\sqrt{b_{j}}.
\end{equation}
If $\sqrt{b_{1}}=\sqrt{b_{3}}$, then $d_{1}=d_{3}$. Multiplying $4a_{j}$ on both sides of (\ref{18}), then $(4\alpha+1)a_{j}\sqrt{b_{j}}\approx 4d_{j}a_{j}$ implies $\Delta_{\lambda_{j}}\approx\frac{4d_{j}a_{j}}{(4\alpha+1)}$. Take $t_{0}=(4\alpha+1)\pi$. By a simple calculation,
$$\cos \frac{{{\Delta _{{\lambda _j}}}t_{0}}}{2}\approx \cos\dfrac{\frac{4d_{j}a_{j}}{(4\alpha+1)}(4\alpha+1)\pi}{2}=\cos2a_{j}d_{j}\pi=1. $$
Thus,
\begin{equation}
\begin{split}
 &\left|{e_{(v_{1},0)}}\exp ( - it_{0}A({C_{4}} \bigstar {G_{2}})){e_{(v_{3},0)}}\right|\\
 &\;\;\;\;\;\;\;\;= \left|\sum\limits_{\lambda_{j}\neq0} {{e^{ - it_{0}\frac{{{\lambda _j} + k}}{2}}}} {e_{v_{1}}}{f_{{\lambda _j}}}({C_4}){e_{v_{3}}}(\cos \frac{{{\Delta _{{\lambda _j}}}t_{0}}}{2} + i\frac{{k - {\lambda _j}}}{{{\Delta _{{\lambda _j}}}}}\sin \frac{{{\Delta _{{\lambda _j}}}t_{0}}}{2}) + {e_{v_{1}}}{f_0}({C_4}){e_{v_{3}}}\right|\\
 &\;\;\;\;\;\;\;\;\approx\left|\sum\limits_{\lambda_{j}\neq0} e^{ - it_{0}\frac{k}{2}}e^{ - it_{0}\frac{\lambda _j}{2}}
 {e_{v_{1}}}{f_{{\lambda _j}}}({C_4}){e_{v_{3}}}+{e_{v_{1}}}{f_0}({C_4}){e_{v_{3}}}\right|.
\end{split}
\end{equation}
Since $e^{ - it_{0}\frac{k}{2}}=1$ whenever $k=0(\text{mod} 4)$, $e^{ - it_{0}\frac{\lambda _j}{2}}=e^{-i\frac{\pi}{2}\lambda_{j}}$ and $e^{-i\frac{\pi}{2}0}=1$, then
\begin{equation}\small
\begin{split}
\left|\sum\limits_{\lambda_{j}\neq0} e^{ - it_{0}\frac{k}{2}}e^{ - it_{0}\frac{\lambda _j}{2}}
 {e_{v_{1}}}{f_{{\lambda _j}}}({C_4}){e_{v_{3}}}+{e_{v_{1}}}{f_0}({C_4}){e_{v_{3}}}\right|&=\left|\sum\limits_{\lambda_{j}\neq0} e^{-i\frac{\pi}{2}\lambda_{j}}
 {e_{v_{1}}}{f_{{\lambda _j}}}({C_4}){e_{v_{3}}}+{e^{-i\frac{\pi}{2}0}e_{v_{1}}}{f_0}({C_4}){e_{v_{3}}}\right|\\
&=\left|\sum\limits_{\lambda_{j}\in S} e^{-i\frac{\pi}{2}\lambda_{j}}
 {e_{v_{1}}}{f_{{\lambda _j}}}({C_4}){e_{v_{3}}}\right|\\
&=1,
\end{split}
\end{equation}
where the last equality holds as $C_{4}$ admits PST at time $\frac{\pi}{2}$ between vertices $v_{1}$ and  $v_{3}$. $\Box$\\
\\
\textbf{Example.} The complete graph $K_{5}$ is a $4$-regular connected graph, then $k=4$ and $n_{2}=5$. By a simple calculation, $\sqrt{(2+k)^{2}+16n_{2}}=\sqrt{116}$ and $\sqrt{(2-k)^{2}+16n_{2}}=\sqrt{84}$. Neither of $\sqrt{(2+k)^{2}+16n_{2}}$ and $\sqrt{(2-k)^{2}+16n_{2}}$ is an integer, then ${C_4} \bigstar {K_{5}}$ admits PGST.

\end{document}